\documentclass[12pt]{article} 

\usepackage[left= 0.8in, textwidth=6.9in, textheight=9.0in, top=1.0in]{geometry}
\usepackage{graphicx}
\usepackage{amsmath,amssymb,verbatim}
\usepackage{enumitem}
\usepackage{color}  
\usepackage{bm}                                  

\newcommand{\MT}{\left[ \begin{array}{rrrrrrrrrrrrrrrrrrrr}}
\newcommand{\EM}{\end{array}\right]}
\newcommand{\EQ}{\begin{equation}\begin{array}{lllllllllll}}
\newcommand{\EE}{\end{array}\end{equation}}
\newcommand{\Real}{\mathbb R}
\newcommand\norm[1]{\left\lVert#1\right\rVert}
\newcommand\abs[1]{\left|#1\right|}

\newtheorem{theorem}{Theorem}
\newtheorem{lemma}{Lemma}

\newtheorem{definition}{Definition}

\newtheorem{corollary}{Corollary}

\def\Fr{\ds \frac}
\def\ds{\displaystyle}
\def\calD{{\cal D}}

\def\calV{{\cal V}}
\def\calE{{\cal E}}

\def\calS{{\cal S}}
\def\calR{{\cal R}}

\def\bff{{\bf f}}

\def\bfx{{\bf x}}
\def\bfz{{\bf z}}

\def\bfw{{\bf w}}
\def\bfW{{\bf W}}

\def\bfb{{\bf b}}

\def\bfphi{{\boldsymbol\phi}}

\title{\LARGE \bf
Data-Driven Computational Methods for the Domain of Attraction and Zubov's Equation
\thanks{This work was supported in part by U.S. Naval Research Laboratory - Monterey, CA}
}

\author{Wei Kang\thanks{Department of Applied Mathematics, Naval Postgraduate School, Monterey, CA, USA; wkang@nps.edu} \thanks{ Department of Applied Mathematics, University of California at Santa Cruz, Santa Cruz, CA, USA}
\and Kai Sun\thanks{Department of EECS, University of Tennessee, Knoxville, TN, USA; kaisun@utk.edu}
\and Liang Xu\thanks{Marine Meteorology Division, Naval Research Laboratory, Monterey, CA, USA; liang.xu@nrlmry.navy.mil}
}

\begin{document}
\maketitle

\abstract{This paper deals with a special type of Lyapunov functions, namely the solution of Zubov's equation. Such a function can be used to characterize the domain of attraction for systems of ordinary differential equations. In Theorem \ref{thm1}, we derive and prove an integral form solution to Zubov's equation. For numerical  computation, we develop two data-driven methods. One is based on the integration of an augmented system of differential equations; and the other one is based on deep learning. The former is effective for systems with a relatively low state space dimension and the latter is developed for high dimensional problems. The deep learning method is applied to a New England 10-generator power system model. A feedforward neural network is trained to approximate the corresponding Zubov’s equation solution. The network characterizes the system's domain of attraction. We prove that a neural network approximation exists for the Lyapunov function of power systems such that the approximation error is a cubic polynomial of the number of generators. The error convergence rate is $O(n^{-1/2})$, where $n$ is the number of neurons. }

\section{Introduction}
For systems of ordinary differential equations (ODEs), it is well known that their stability can be determined by using Lyapunov functions. In addition, a level set of a Lyapunov function forms the boundary of an invariant set inside the domain of attraction (DOA). Originally introduced in the doctoral dissertation of Lyapunov, the theory has been actively studied for more than a century with applications to almost all fields of science and engineering that deal with dynamic phenomena. There have existed tremendous works of literature on the study of how to find Lyapunov functions for systems of ODEs. Some widely used methodologies are based on the results and tools in related areas such as Zubov's equation \cite{ferguson,giesl,gruene,kinnen,vannelli,zubov},  linear programming (LP) \cite{johansen,julian} and linear matrix inequality (LMI) \cite{bartels,boyd,chesi,mikkelsen,parrilo,polik}. In these works, the Lyapunov functions are represented using various types of approximations, including quadratic or high order polynomials \cite{bellman,ferguson,zubov}, rational functions \cite{vannelli}, collocation \cite{giesl}, piecewise continuous functions \cite{goebel,julian} and sum-of-squares (SOS) \cite{ahmadi,anderson,parrilo}. A comprehensive literature review is beyond the scope of this paper. For the interested readers, many survey papers were published in the last few decades summarizing research advancements at different phases in the development of Lyapunov function theory and its computation \cite{anderson,chesi,chen,giesl1,gurel,polik}. 

Existing methods of solving Zubov's equation have some fundamental limitations. Notably, most methods suffer the curse of dimensionality, i.e., the complexity of the solution increases exponentially with the dimension of the state space. In this case, the number of grid points in discretization, or coefficients in polynomial approximation, increases so fast that the computational algorithm is intractable in high dimensional state spaces. For Lyapunov functions, techniques of using special forms of representation, such as piecewise quadratic function, rational function and SOS, cannot characterize the true boundary of the DOA although they can be used to find invariant subsets inside the DOA.  It is worth noting that ideas of using deep learning to mitigate the curse of dimensionality have been successfully applied to various PDEs over the past few years \cite{han,nakagongkang1,raissi}. Its potential to applied problems that have high dimensions, such as numerical weather prediction,  has attracted increasing attention \cite{kochkov,otsuka,shi}. For the specific problem of Lyapunov function and Zubov's equation, a method based on deep learning is studied in \cite{gruene} under the assumption that the system satisfies a small-gain condition and it admits a compositional Lyapunov function. 
 
In this paper, we introduce a new method of solving Zubov's equation, a partial differential equation (PDE) associated with dynamical systems. The solution of this PDE  is a special Lyapunov function that can be used to characterize the DOA as well as its boundary. The contributions of this paper include: (a) deriving and proving an integral form solution to Zubov's equation for general systems of ODEs; (b) developing two data-driven computational methods; (c) as a part of an applied example, proving a theorem about the error upper bound and neural network complexity for a power system model. One of the computational methods developed in this paper is based on integrating an augmented system of ODEs, and the other one is based on deep learning. The former computes the Lyapunov function using ODE solvers, which is effective for relatively low dimensional problems.  We introduce deep neural networks to approximate the Lyapunov function if a system has a high dimension. For example, the computational methods are applied to a 10-generator power system, which is a reduced model of the power grid in New England. A feedforward neural network is trained that characterizes the DOA of the power system, which is an interesting research result in its own right.

\section{An integral form solution to Zubov's equation}
\label{sec_zubov}
In the stability theory for dynamical systems, it is known that a solution of Zubov's equation is a special kind of Lyapunov function. Given any point in the state space, the point is located inside the DOA if and only if the value of the Lyapunov function is less than one. In this section, we prove that a function defined in an integral form solves Zubov's equation. Consider the following dynamical system, 
\EQ
\label{eq_sys}
\dot \bfx (t)= \bff (\bfx (t)), & \bfx \in \Real^n
\EE
where $\bfx$ is the state variable and $n$ is the dimension of the state space. The solution of (\ref{eq_sys}) is denoted by 
\EQ
\bfx (t)=\bfphi (t, \bfx_0),
\EE
where $\bfx_0$ is the initial state, $\bfphi(0,\bfx_0)=\bfx_0$. We assume that, for any initial state, the solution of (\ref{eq_sys}) is unique and defined in $t\in [0, \infty)$.  A point $\bfx\in \Real^n$ is said to be an equilibrium if 
$\bff(\bfx)=0$.
An invariant set, denoted by $\calE$, is a subset of $\Real^n$ so that
\EQ
\bfphi(t,\bfx)\in \calE, & t\geq 0,
\EE
 for arbitrary $\bfx\in \calE$. The concept of invariant set is a generalization of equilibrium because a set consisting of equilibrium points is an invariant set.  

In this paper, $\norm{\bfx}$ represents the Euclidean norm on $\Real^n$. The distance between a point and a set is defined as the infimum of distances between the point and all points in the set, denoted by $d(\cdot, \cdot)$. For example, the distance from a point, $\bfx_0$, to an invariant set, $\calE$, is
\EQ
d(\bfx_0, \calE)=\inf\{\norm{\bfx_0-\bfx}; \; \bfx\in \calE\}.
\EE
Given a positive number $\delta$,  the $\delta$-neighborhood of $\calE$ is defined by
\EQ
B_\delta (\calE)=\{ \bfx\in\Real^n; \; d(\bfx, \calE) <\delta \}.
\EE

\begin{definition}
(Asymptotic stability) \\
(a) A closed invariant set $\calE\subset \Real^n$ of (\ref{eq_sys}) is called asymptotically stable if the following assumptions hold true. 
\begin{enumerate}
\item For any $\epsilon > 0$, there exists a $\delta >0$ such that $\bfx\in B_\delta(\calE)$ implies $d(\bfphi(t,\bfx),\calE)<\epsilon$ for all $t\geq 0$.
\item There exists a $\delta >0$ such that
\EQ
\label{eq_asymp}
 \ds\lim_{t\rightarrow \infty} d(\bfphi(t,\bfx),\calE) = 0.
\EE
whenever $\bfx\in B_\delta(\calE)$.
\end{enumerate} 
The domain of attraction is denoted by $\calD$,
\EQ
\label{eq_D}
\calD = \{ \bfx \in \Real^n; \; \ds\lim_{t\rightarrow \infty} d(\bfphi(t,\bfx),\calE) = 0\}.
\EE
It is an open set \cite{hahn}. The closure of $\calD$ is denoted by $\bar\calD$.\\
(b) The set $\calE$ is called uniformly asymptotically stable if, in addition, there exists a $\delta$-neighborhood of $\calE$ such that the limit (\ref{eq_asymp}) converges uniformly on $B_\delta(\calE)$. More specifically, for any $\epsilon>0$, there exists $T>0$ such that 
\EQ
d(\bfphi(t,\bfx),\calE)<\epsilon
\EE
for all  $\bfx\in B_\delta (\calE)$ and $t\geq T$. \\
(c) The set $\calE$ is called uniformly attracting if there exits $\delta>0$ such that, for any $0<h<\delta$, there exist $T>0$ and $\gamma>0$. Whenever $h<d(\bfx,\calE)<\delta$ and $t\in [0, T]$, we have 
\EQ
d(\bfphi(t,\bfx),\calE)>\gamma.
\EE 
\end{definition}

\noindent \textbf{Remark} For autonomous systems ($\bff(\bfx)$ is independent of $t$), a compact asymptotically stable $\calE$  must be uniformly asymptotically stable \cite{zubov}. Or if $\bfphi(t,\bfx)$ converges to $\calE$ exponentially for $\bfx\in B_\delta (\calE)$, i.e.,
\EQ
\label{eq_uniformstable}
d(\bfphi (t,\bfx),\calE)\leq M d(\bfx,\calE) e^{-\lambda t}
\EE
for some $M>0$ and $\lambda >0$, then $\calE$ is uniformly asymptotically stable. For problems considered in this paper, it is justified in the proof of Theorem \ref{thm1} that $\calE$ is always uniformly attracting. \\

\noindent \textbf{Assumption 1}. We assume that $\bff(\bfx)$ in (\ref{eq_sys}) is locally Lipschitz on $\Real^n$ and $\bfphi(t,\bfx)$ exists for all $(t,\bfx)\in[0, \infty)\times \Real^n$. Furthermore, $\norm{\bff(\bfx)}$ is bounded on $B_\delta (\calE)$ for some $\delta >0$. \\

\noindent \textbf{Assumption 2} $\,$Let $\calE$ be a closed invariant set of (\ref{eq_sys}). Let $W: \Real^n\rightarrow \Real$ be a continuous function. We assume 
\EQ
\label{eq_W}
W(\bfx)=0,& \mbox{if } \bfx\in\calE,\\ 
W(\bfx)> 0, & \mbox{otherwise},
\EE
and
\EQ
\label{eq_Wlimit}
\ds\lim_{d(\bfx,\calE)\rightarrow 0}W(\bfx)=0.
\EE
In addition, for any $\delta>0$, there exists a $\gamma>0$ such that
\EQ
\label{eq_Wbound}
W(\bfx)>\gamma
\EE
whenever $d(\bfx,\calE)> \delta$.\\

If $\calE$ is a compact set, then (\ref{eq_Wlimit}) is implied by (\ref{eq_W}) and the continuity of $W(\bfx)$. Among the functions satisfying Assumption 2, two useful ones that we adopt for the examples in this paper are

\EQ
&W(\bfx)=\norm{\bfx -\bfx_e}^2,  & \mbox{ if } \calE=\{ x_e\}, \mbox{ where } \bfx_e \mbox{ is a known equilibrium,}\\
\mbox{or} &W(\bfx)=\norm{\bff(\bfx)}^2, & \mbox{ if the system has multiple equilibrium points}
\EE
The former is a simple function that can be used for problems with a known equilibrium. The latter is applicable to not only isolated equilibrium but also systems that have a set of equilibrium points. The following integral function is used for the construction of a solution to Zubov's equation,
\EQ
\label{eq_Ix}
I(\bfx)= \ds\int_0^\infty W(\bfphi (s,\bfx))ds
\EE
where $\bfphi(t,\bfx)$ is the solution of (\ref{eq_sys}) satisfying $\bfphi(0,\bfx)=\bfx$. \\

\noindent \textbf{Assumption 3} $\,$Let $\calE$ be a closed invariant set. We assume that the integral in (\ref{eq_Ix}) converges in a neighborhood $\bfx\in B_\delta(\calE)$ for some $\delta>0$. Furthermore, we assume
\EQ
\label{eq_Ixlimit}
\ds\lim_{\bfx\in B_\delta(\calE),d(\bfx,\calE)\rightarrow 0}I(\bfx)=0.
\EE

\noindent \textbf{Remark} $\;$ If $W(\bfx)$ is Lipschitz in $B_\delta (\calE)$ for some $\delta >0$, if $\bfphi(t,\bfx)$ converges to $\calE$ exponentially in the sense of (\ref{eq_uniformstable}), then Assumption 3 holds. To prove this claim, let $L$ be the Lipschitz constant of $W(\bfx)$, then
\EQ
W(\bfphi(t,\bfx)) \leq L \norm{\bfphi(t,\bfx) -\bar\bfx}
\EE
for all $\bar\bfx\in \calE$. Therefore,
\EQ
I(\bfx)=\ds\int_0^\infty W(\bfphi(s,\bfx))ds 
\leq \ds\int_0^\infty  L d(\bfphi(t,\bfx),\calE)ds \leq \ds\int_0^\infty  M d(\bfx,\calE) e^{-\lambda s}ds=\Fr{M}{\lambda}d(\bfx,\calE)
\EE
Therefore, $I(\bfx)$ is a finite number; and it satisfies (\ref{eq_Ixlimit}). \\

Zubov's method is summarized in the following theorem. Its proof can be found in \cite{zubov}.

\begin{theorem}
\label{thm_zubov}
(Zubov's Theorem \cite{zubov}) Let $\calE$ be a closed invariant set of (\ref{eq_sys}) in which $\bff(\bfx)$ satisfies Assumption 1.  Suppose $\calE$ is uniformly asymptotically stable and uniformly attracting. Let $\calD$ be an open set that contains a $\delta$-neighborhood of $\calE$. Then $\calD$ is the domain of attraction if and only if there exist continuous functions $V(\bfx)$ defined on $\calD$ and $\Psi(\bfx)$  defined on $\Real^n$ that satisfy the following properties:
\begin{enumerate}[label={(\arabic*)}]
\item If $\bfx\in \calD$ and $d(\bfx, \calE) \neq 0$, then $0<V(\bfx)<1$. 
\item If $\bfx\in \Real^n$ and $d(\bfx, \calE) \neq 0$, then  $\Psi(\bfx)>0$. If $\bfx \in \calE$, then $\Psi (\bfx)=0$. 
\item If $d(\bfx,\calE)\rightarrow 0$, then $V(\bfx)\rightarrow 0$ and $\Psi(\bfx)\rightarrow 0$.
\item For any sufficiently small $\delta>0$, there exist $\gamma_1>0$ and $\gamma_2>0$ such that $d(\bfx,\calE)\geq \delta$ implies
\EQ
\label{eq_cond3}
V(\bfx)>\gamma_1 \mbox{ and } 
\Psi(\bfx) > \gamma_2
\EE
\item For any $\bar \bfx\in \bar \calD\setminus \calD$, 
\EQ
\ds\lim_{\bfx\in \calD, \bfx\rightarrow \bar \bfx} V(\bfx)=1
\EE
\item For any $\bfx\in\calD$, $V(\bfphi(t,\bfx))$ is differentiable with respect to $t$ satisfying,
\EQ
\left( \Fr{dV(\bfphi(t,\bfx))}{dt}\right)_{t=0} = -\Psi(\bfx)(1-V(\bfx)).
\EE
\end{enumerate}
\end{theorem}

In the rest of this section, we will prove that the following function forms a solution to Zubov's equation,
\EQ
\label{eq_soln_zubuv}
V(\bfx)=\left\{ \begin{array}{lll}
\tanh (\alpha I(\bfx)), &\mbox{ if } I(\bfx)<\infty,\\
1, & \mbox{ otherwise},
\end{array}\right.\\
\Psi (\bfx)=\alpha(1+V(\bfx))W(\bfx),
\EE
where $\alpha >0$ is a constant. 
\begin{theorem}
\label{thm1}
Consider a system (\ref{eq_sys}) and a closed invariant set $\calE$. Suppose $\calE$ is uniformly asymptotically stable. Suppose $\bff(\bfx)$ satisfies Assumption 1. Suppose $W: \Real^n \rightarrow \Real$ is a function satisfying Assumptions 2 and 3. Then, $V(\bfx)$ and $\Psi(\bfx)$ defined in (\ref{eq_soln_zubuv}) solves Zubov's equation in the sense that they satisfy items (1)-(6) in Theorem \ref{thm_zubov}.
\end{theorem}

Some properties of $I(\bfx)$ is essential. We prove the following lemma first, before the proof of Theorem \ref{thm1}. 

\begin{lemma} 
\label{lemma1}
Under the same assumption as in Theorem \ref{thm1}, $I(\bfx)$ has the following properties. 
\begin{enumerate}[label={(\arabic*)}]
\item $I(\bfx) < \infty$ for $\bfx\in\calD$ and $I(\bfx)$ is a continuous function on $\calD$. 
\item $I(\bfx) >0$ if $\bfx\in \calD\setminus \calE$ and $I(\bfx) =0$ if $\bfx\in\calE$.
\item If $\bfx\notin \calD$, $I(\bfx)=\infty$. If $\bar \bfx\in \bar \calD\setminus \calD$, then 
\EQ
\label{eq_claim3}
\ds\lim_{\bfx\in \calD, \bfx\rightarrow \bar\bfx}I(\bfx)=\infty
\EE
\item Given any $\bfx\in \calD$, 
\EQ
\Fr{d\,}{dt}I(\bfphi(t,\bfx)) = -W(\bfphi (t,\bfx))
\EE
\end{enumerate}
\end{lemma}

\noindent{\it Proof}. 

(1) For any $\bfx \in \calD$ and a radius, $\delta$, defined in Assumption 3, (\ref{eq_D}) implies that there exits a $t_1>0$ such that 
\EQ
d(\bfphi (t,\bfx),\calE) < \delta
\EE
if $t\geq t_1$. From Assumption 3, the integral
$$\ds\int_0^\infty W(\bfphi (s,\bfphi(t_1,\bfx)))ds, $$
converges. Therefore, 
\EQ
I(\bfx)=\ds\int_0^\infty W(\bfphi (s,\bfx))ds=\ds\int_0^{t_1} W(\bfphi (s,\bfx))ds+\ds\int_0^\infty W(\bfphi (s,\bfphi(t_1,\bfx)))ds < \infty.
\EE
In the next, we prove that $I(\bfx)$ is continuous on $\calD$. Define
\EQ
\label{eq_ztx}
z(t, \bfx)=\ds\int_0^t W(\bfphi (s,\bfx))ds
\EE
Because $\bfphi(t,\bfx)$ is continuous, it is a known fact in real analysis that $z(t,\bfx)$ is a continuous function of $\bfx$ for any fixed value of $t$. For any $\epsilon>0$, Assumption 3 implies that  
\EQ
\label{eq_intT}
I(\bfx)=\ds\int_0^\infty W(\bfphi (s,\bfx))ds <\epsilon, 
\EE
if $\bfx \in B_{\delta_1}(\calE)$ for some $\delta_1>0$. Given any $\bfx\in\calD$, there exists a $t_1>0$ such that $\bfphi(t_1,\bfx)\in B_{\delta_1}(\calE)$. Because $\bfphi(t_1,\bfx)$ is a continuous function of $\bfx$, there exists $\delta_2 >0$ such that $\bfphi(t_1,\bfz)\in B_{\delta_1}(\calE)$ if $\norm{\bfz-\bfx}<\delta_2$. 
From (\ref{eq_intT}) and the fact that $z(t,\bfx)$ is a continuous function of $\bfx$, we have
\EQ
\label{eq_Ixcont}
\ds\lim_{\bfz\rightarrow \bfx}\abs{I(\bfz)-I(\bfx)}\\
=\ds\lim_{\bfz\rightarrow \bfx}\abs{z(t_1,\bfz)+\ds\int_{0}^\infty W(\bfphi (s,\bfphi(t_1,\bfz)))ds-z(t_1,\bfx)-\ds\int_{0}^\infty W(\bfphi (s,\bfphi(t_1,\bfx)))ds}\\
\leq \ds\lim_{\bfz\rightarrow \bfx}\abs{z(t_1,\bfz)-z(t_1,\bfx)}+2\epsilon\\
=2\epsilon
\EE
Because $\epsilon$ can be arbitrarily small, (\ref{eq_Ixcont}) implies that $I(\bfx)$ is continuous at every $\bfx\in \calD$.

(2) Let $\bfx$ be a point in $\calD\setminus \calE$. From (\ref{eq_W}), $W(\bfphi(t,\bfx))>0$ for small value of $t>0$. Therefore, its integral is greater than zero, i.e., $I(\bfx)>0$. If $\bfx\in\calE$, we have $\bfphi(t,\bfx)\in \calE$ for $t > 0$ because $\calE$ is invariant.  Then (\ref{eq_W}) implies $W(\bfphi(t,\bfx))\equiv 0$ and $I(\bfx)=0$.

(3)  If $\bfx\notin \calD$, then $\bfphi(t,\bfx)\notin \calD$ for all $t \geq 0$ (otherwise, $\bfphi(t,\bfx)$ approaches $\calE$). From Assumption 2 and the assumption about $\calD$ in Theorem \ref{thm_zubov}, there exist numbers $\delta>0$ and $\gamma>0$ such that $B_\delta(\calE)\subset \calD$ and
\EQ
W(\bfx)>\gamma
\EE
whenever $d(\bfx,\calE)> \delta$. Therefore, 
\EQ
W(\bfphi(t,\bfx))> \gamma
\EE
for all $t\geq 0$. Therefore, 
\EQ
I(\bfx)=\ds\int_0^\infty W(\bfphi(s,\bfx))ds =\infty
\EE
For any $\bar\bfx \in \bar\calD \setminus \calD$, $\bar\bfx \notin \calD$. Therefore, $I(\bar\bfx)=\infty$, i.e., for any $M>0$, there exists a $T>0$ such that 
\EQ
\label{eq_int0T}
\ds\int_0^T W(\bfphi(s,\bar\bfx))ds > 2M
\EE
It is known that the integral in (\ref{eq_int0T}) is continuously dependent on $\bar\bfx$. There exists an open neighborhood of $\bar \bfx$ such that, for any $\bfx$ in the neighborhood,  
\EQ
\label{eq_lim_M}
\ds\int_0^T W(\bfphi(s,\bfx))ds > M
\EE
Therefore, for any $M>0$, there exists an open neighborhood of $\bar\bfx$ in which 
\EQ
I(\bfx)\geq \ds\int_0^T W(\bfphi(s,\bfx))ds>M,
\EE
i.e., (\ref{eq_claim3}) is proved. 

(4) Given any $\bfx\in \calD$, we have
\EQ
\Fr{d\,}{dt}I(\bfphi(t,\bfx)) \\
= \Fr{d\,}{dt} \left( \ds\int_0^\infty W(\bfphi (s,\bfphi(t,\bfx)))ds \right)\\
=\Fr{d\,}{dt} \left( \ds\int_t^\infty W(\bfphi (s,\bfx))ds \right)\\
=-W(\bfphi (t,\bfx)).
\EE
$\blacksquare$\\

\noindent {\it Proof of Theorem \ref{thm1}}. In the following, we prove that $V(\bfx)$ and $\Psi(\bfx)$ satisfy items (1)-(6) in Theorem \ref{thm_zubov}. 
For general dynamical systems, it is required in Theorem \ref{thm_zubov} that $\calE$ is uniformly attracting. In this paper, all dynamical systems are defined by equations in the form of (\ref{eq_sys}) in which $\bff(\bfx)$ is continuous. If Assumption 1 holds, then $\norm{\bff(\bfx)}$ is bounded on  $0<h \leq d(\bfx,\calE) \leq \delta$. If $T>0$ is small enough, then $d(\bfphi(t,\bfx),\calE)$  is uniformly bounded from below by a $\gamma >0$. Therefore, $\calE$ is uniformly attracting. 

From (\ref{eq_Ix}), $V(\bfx)\equiv 1$ is continuous in the interior of $\Real^n\setminus\calD$. 
From Lemma \ref{lemma1}(1), $V(\bfx)$ is continuous on $\calD$. From  Lemma \ref{lemma1}(3), $V(\bfx)$ is continuous at $\bfx$ on the boundary of $\calD$. Therefore, $V(\bfx)$ is continuous on $\Real^n$. Consequently, $\Psi (\bfx)$ is continuous on $\Real^n$. 

Lemma \ref{lemma1}(1)-(2) as well as the facts that $\tanh(0)=0$ and $0<\tanh(x)<1$ for $0< x<\infty$ imply $0< V(\bfx) <1$ if $\bfx\in\calD\setminus\calE$  and $V(\bfx)=0$ if $\bfx\in\calE$. Therefore, $V(\bfx)$ satisfy Condition (1) in Theorem \ref{thm_zubov}. 

If $\bfx\notin\calE$, then $W(\bfx)> 0$. Therefore, $\Psi(\bfx)>0$. If $\bfx\in \calE$, then $W(\bfx)=0$. Therefore, $\Psi(\bfx)=0$. Condition (2) in Theorem \ref{thm_zubov} is satisfied. 

Because of (\ref{eq_Wlimit}) and (\ref{eq_Ixlimit}), we know $\Psi(\bfx)\rightarrow 0$ and $V(\bfx)\rightarrow 0$ as $d(\bfx,\calE)\rightarrow 0$. 
Therefore, $V(\bfx)$ and $\Psi(\bfx)$ satisfy Condition (3) in Theorem \ref{thm_zubov}. 

From (\ref{eq_Wbound}) in Assumption 2, $\Psi(\bfx)$ satisfies Condition (4). From Assumption 1, $\norm{\bff(\bfx)}$ is bounded on $B_\delta(\calE)$ for any $\delta >0$ that is small enough. Therefore, there is a $T>0$ such that $d(\bfphi(T,\bfx), \calE)>\delta/2$ if $d(\bfx,\calE)>\delta$. From (\ref{eq_Wbound}), we know 
\EQ
I(\bfx)\geq \ds\int_0^T W(\bfphi(s,\bfx))ds \geq \gamma T
\EE
for some $\gamma >0$, provided $d(\bfx,\calE)>\delta$. Therefore, $V(\bfx)$ satisfies Condition (4). 

Lemma \ref{lemma1}(3) implies Condition (5). 

Now, we can prove that $V(\bfx)$ and $\Psi(\bfx)$ solves Zubov's equation,
\EQ
 \left( \Fr{dV(\bfphi(t,\bfx))}{dt}\right)_{t=0}=\left. (1-\tanh^2(\alpha I(\bfphi(t,\bfx))))\alpha\left( \Fr{dI(\bfphi(t,\bfx))}{dt}\right)\right|_{t=0}\\
 =\left. -\alpha\left(1-V^2(\bfphi(t,\bfx)))W(\bfphi(t,\bfx)\right)\right|_{t=0}\\
 =-\alpha (1+V(\bfx))W(\bfx)(1-V(\bfx))\\
 =-\Psi(\bfx)(1-V(\bfx))
\EE 
$\blacksquare$

\section{A computational algorithm based on ODE solvers}
\label{sec_algorithm}

The function, $V(\bfx)$, is a global Lyapunov function which defines the domain of attraction. 
\begin{corollary}
Under the same assumption as in Theorem \ref{thm1}, 
$\bfx\in \calD$ if and only if $V(\bfx) <1$. 
\end{corollary}

\begin{figure}[!ht]
\centering
\includegraphics[width = 4.0in]{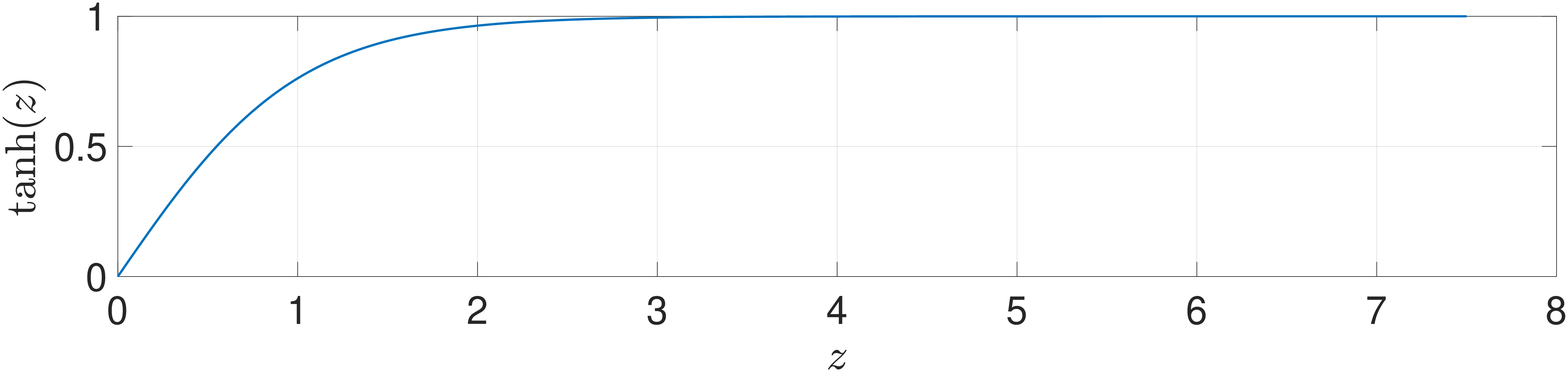}
\caption{ $\tanh(z)$ }
\label{fig_tanh}
\end{figure}

To characterize the DOA in this way, we must find effective and efficient computational algorithms for $V(\bfx)$. It is known from (\ref{eq_claim3}) that, when $\bfx$ approaches the boundary of $\calD$, $I(\bfx)\rightarrow\infty$. The computation is increasingly time consuming in a thin layer along the boundary of $\calD$. To be more specific, a {\it boundary layer} of $\calD$ is defined as follows 
\EQ
\calD_{layer}(\delta) = \{ \bfx \in \calD; 1-\delta < V(\bfx) < 1 \}
\EE 
where $0<\delta<1$. In (\ref{eq_soln_zubuv}), $V(\bfx)$ is defined based on $\tanh(z)$, which is shown in Figure \ref{fig_tanh}. It is convenient to define the value of $\delta$ using $\tanh(z)$. Because $1-\tanh(20)<1\times 10^{-16}$ is around machine precision, the boundary layer $\calD_{layer}(\tanh(20))$ is numerically indistinguishable from the boundary of $\calD$. This boundary layer is used to determine the scale factor, $\alpha$, in (\ref{eq_soln_zubuv}). Details are illustrated in Example 1.

To approximate $I(\bfx)$ for $\bfx$ in $\calD$, the integral in (\ref{eq_Ix}) is approximated by integrating over a finite interval $[0, T]$, where $T$ is large enough so that 
\EQ
\left. \Fr{d\;}{dt}\ds\int_0^t W(\bfphi (s,\bfx))ds \right|_{t=T}< \delta_I
\EE
where $\delta_I>0$ is a small number. For $\bfx$ located outside the DOA, $I(\bfx)=\infty$. Numerically, we need a number $M>0$ that is large enough to terminate the integration process, i.e. if 
\EQ
\label{eq_largeM}
\ds\int_0^T W(\bfphi (s,\bfx))ds > M
\EE
we terminate the integration. In this case, the value of $I(\bfx)$ is considered to be $\infty$. Note that (\ref{eq_largeM}) does not guarantee that $\bfx$ is outside the DOA, no matter how large is the value of $M$. However, if the value of $M$ is large enough, the set inside the DOA satisfying $I(\bfx)>M$ is contained in a thin boundary layer. In this case, the trajectory is practically unstable (a small perturbation may cause the trajectory to diverge). The following system of equations is called the augmented system of (\ref{eq_sys}),
\EQ
\label{eq_sys_aug}
\dot \bfx (t)= \bff (\bfx (t)), & \bfx(0) \in \Real^n,\\
\dot z=W(\bfx), & z(0)=0.
\EE
Note that the initial condition of the augmented system is in $\Real^n$, not $\Real^{n+1}$, because $z(0)=0$ is a fixed initial value. In the following computational algorithm, $\delta_I$, $M$ and $\alpha$ are positive constants in which  $\delta_I$ is the stopping criteria for the computation of $I(\bfx)$. The values of $M$ and $\alpha$ are determined using a data-driven approach that is illustrated in Example 1.\\

\noindent \textbf{Algorithm 1} $\;$ 
Let $\delta_I$, $M$, $\alpha$ and $\varDelta T$ be positive numbers. Given $\bfx \in\Real^n$, $V(\bfx)$ is computed as follows.
\begin{enumerate}
\item Solving (\ref{eq_sys_aug}) using an ODE solver for $t\in [(k-1)\varDelta T, k\varDelta T]$, $k=1,2,...K$, satisfying $\bfx(0)=\bfx$. The initial state of the $k$th solution is the final state of the $(k-1)$th solution. Stop the process at $K$ so that the $K$th solution satisfies either
\EQ
\label{eq_dzdt1}
z(K\varDelta T) > M
\EE
or
\EQ
\label{eq_dzdt}
\left. \Fr{\varDelta z}{\varDelta t}\right|_{t=K\varDelta T} < \delta_I 
\EE

In (\ref{eq_dzdt}), $\varDelta t $ is the time step of the ODE solver, $\varDelta z=z(K\varDelta T)-z(K\varDelta T -\varDelta t)$.  In this case, $I(\bfx)\approx z(K\varDelta T)$. 
\item 
\EQ
\label{eq_evalV}
V(\bfx)=\left\{
\begin{array}{lll}
1, & \mbox{if (\ref{eq_dzdt1}) holds}, \\
\tanh(\alpha z(K\varDelta T)), &\mbox{if (\ref{eq_dzdt}) holds}.
\end{array}\right.
\EE
\end{enumerate}

Before the end of this section, we introduce the van der Pol equation as an example of implementing the computational algorithm. \\

\begin{figure}[!ht]
\centering
\includegraphics[width = 4.0in]{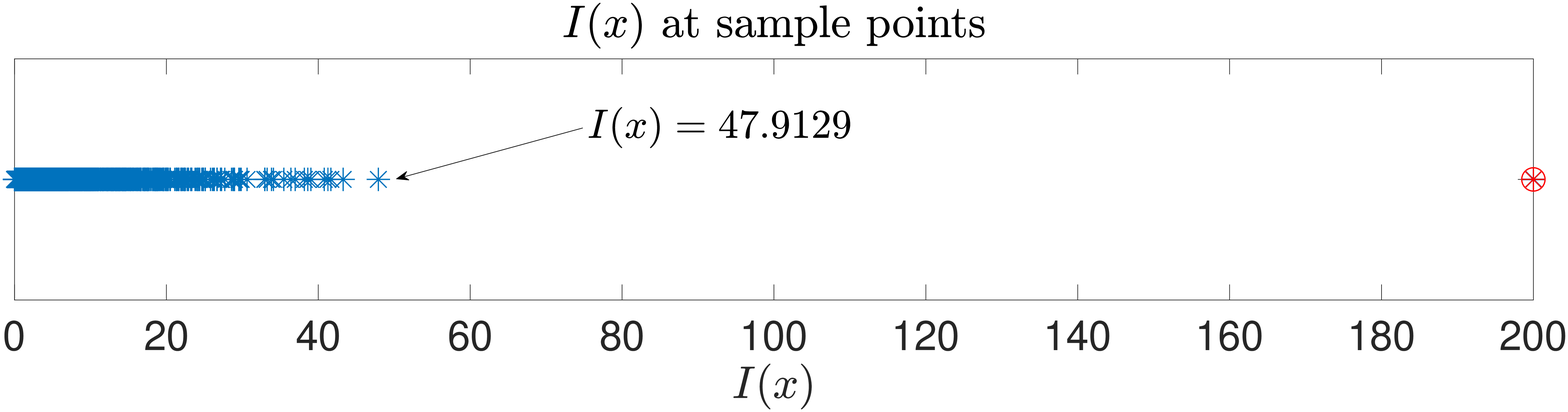}\\
\begin{minipage}{4.5in}
\caption{ I-value plot. The cluster on the left side contains all $I(\bfx_i)$ satisfying (\ref{eq_dzdt}). For sample points satisfying (\ref{eq_dzdt1}), they are symbolically represented by `$ {\scriptstyle \circledast}$' at $I(\bfx)=200$.}
\end{minipage}
\label{fig_Idistribution}
\end{figure}

\begin{figure}[!ht]
\centering
\includegraphics[width = 3.5in]{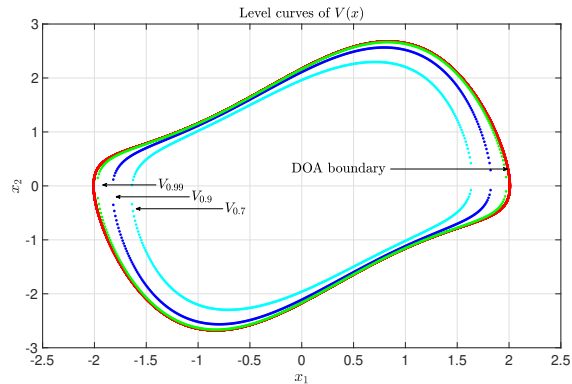}
\caption{ Level curves: $V(\bfx)=r$, $r=0.7$, $0.8$ and $0.99$.}
\label{fig_levelcurves}
\end{figure}

\noindent \textbf{Example 1} $\;$ Consider the following van der Pol equation
\EQ
\dot x_1=-x_2\\
\dot x_2=x_1-(1-x_1^2)x_2
\EE
The origin, $\bfx = 0$, is an asymptotically stable equilibrium. For the computation of $I(\bfx)$, $\delta_{I}=10^{-6}$. A data-driven approach is adopted for the determination of $\alpha$ and $M$. We first generated a set of samples following the uniform distribution, 
\EQ
\calS = \{ \bfx_i; \; 1\leq i\leq N\}
\EE 
where $N=3,000$ is the sample size. Applying Step 1 in Algorithm 1 to the augmented system to compute the value of $I(\bfx)$ for $\bfx\in \calS$. Shown in Figure \ref{fig_Idistribution}, the value of $M=200$ is selected to be large enough so that a gap exists between $M$ and the cluster of $I(\bfx_i)$ satisfying (\ref{eq_dzdt}). We call this graph {\it I-value plot}. In Figure \ref{fig_Idistribution}, all values of $I(\bfx_i)$ satisfying (\ref{eq_dzdt1}) are represented by a single point at $I(\bfx)=200$. It is numerically treated as $\infty$ in the computation of $I(\bfx)$. Increasing the value of $M$ implies the increase of $K\varDelta T$ in (\ref{eq_dzdt1}). The computation will take longer time if $\bfx\notin \calD$. However, if $M$ is too small, some points in $\calD$ is assigned the value $I(\bfx)=\infty$ or $V(\bfx)=1$, i.e., a point inside the DOA is classified as unstable. To avoid that, we choose the value of $M$ that should be larger than $I(\bfx)$ for most states in $\calD$. From the I-value plot, we see that any $M>80$ has a clear gap from the data cluster on the left. By doing that, we reduce the possibility of assigning $V(\bfx)=1$ for $x\in\calD$. In this example, we choose $M=200$. This value should correspond to the boundary layer $\calD_{layer}(\tanh(20))$. For this purpose, we adopt the following scale factor in the computation of $V(\bfx)$
\EQ
\alpha = \Fr{20}{M}.
\EE
Using the value of $\delta_I$, $M$ and $\alpha$ selected as above, $V(\bfx)$ can be computed using Algorithm1. The augmented system is solved using ODE45 in MATLAB. The averaged CPU time of evaluating $V(\bfx)$ using a 2.4 GHz Intel Core i9 processor is $0.03$s. The DOA is characterized by $V(\bfx)<1$. Its boundary can be approximated by a curve $V(\bfx)=r$ if $r<1$ is very close to $1$. Three level curves of $V(\bfx)$ for $r=0.7$, $0.8$ and $0.99$ are shown in Figure \ref{fig_levelcurves}. 

\section{A data-driven method based on deep learning}
The state space of the van der Pol equation in Example 1 is $\Real^2$. Algorithm 1 is effective, although not super fast, for this kind of problems that have a low state space dimension. For ODEs that require significant computational time to solve, typically high dimensional systems,  Algorithm 1 is slow. However, Algorithm 1 has a special property namely {\it causality-free}, i.e., the solution of Zubov's equation is computed at individual points without using the value of $V(\bfx)$ at other points in its neighborhood. Consequently, the computation does not require a grid. Demonstrated in several studies \cite{kang1,kang2,nakagongkang1,nakagongkang2}, causality-free algorithms are advantageous in data-driven computational methods for PDEs. They have perfect parallelism, a desirable property for generating a large amount of data. In contrast to PDE algorithms that are based on grids in space, a causality-free algorithm can solve high dimensional problems. In addition, $V(\bfx)$ can be computed in any targeted area following any distribution that one may choose. The combination of a causality-free algorithm with supervised learning was proved to be effective and efficient for solving high dimensional  PDEs in optimal control \cite{nakagongkang1}.  In this section, we apply supervised learning and feedforward neural networks for the purpose of learning $V(\bfx)$ based on date generated using Algorithm 1. Once the training is done, $V(\bfx)$ is approximated by a neural network. Its computation is much faster than Algorithm 1. \\

\begin{figure}[!ht]
\centering
\includegraphics[width = 3.25in]{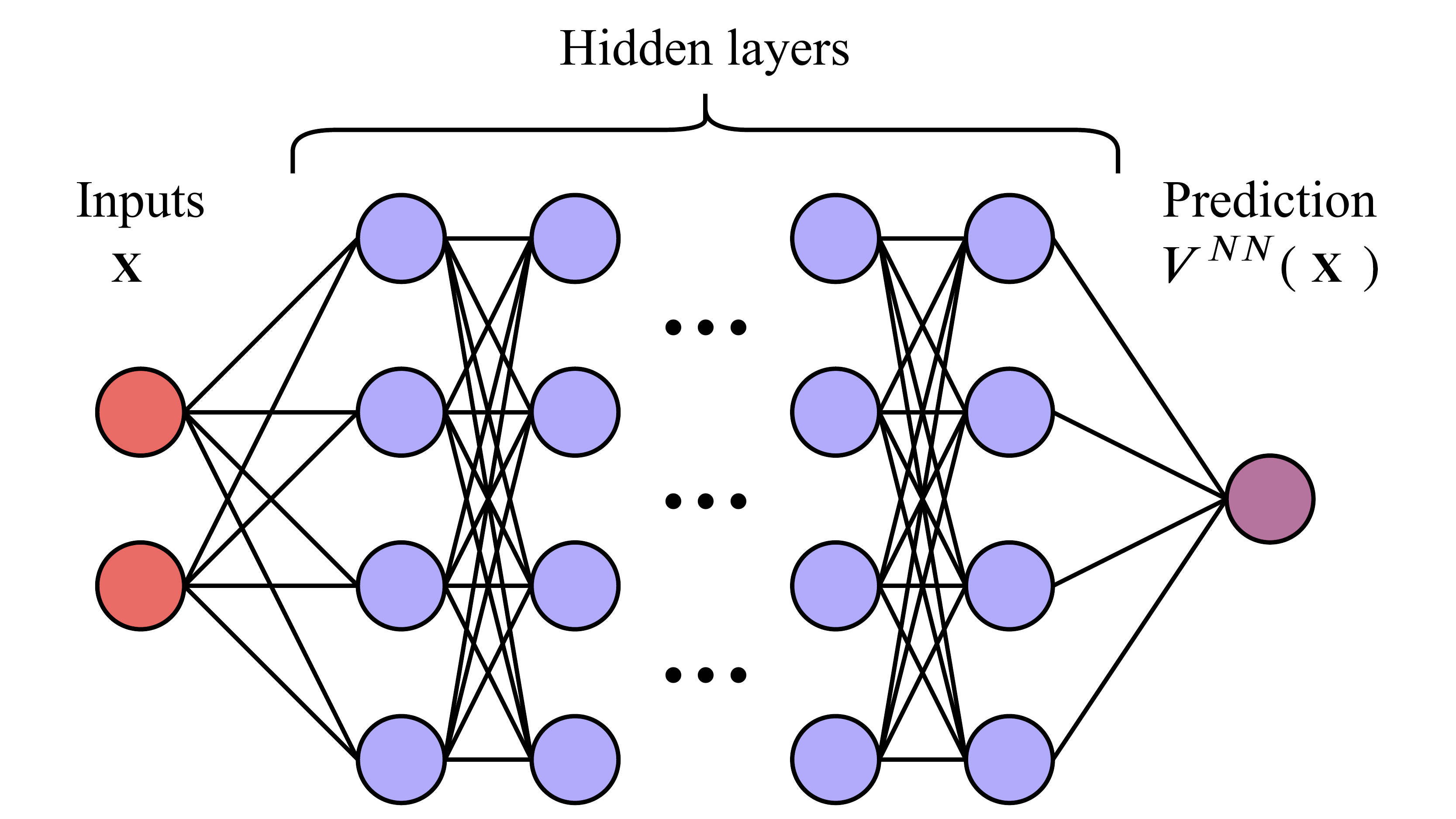}
\caption{ A feedforward neural network }
\label{fig_NNplot}
\end{figure}

\noindent \textbf{Algorithm 2} $\;$ Consider a system (\ref{eq_sys}), a closed invariant set $\calE$, and a function $W: \Real^n \rightarrow \Real$. Suppose that they satisfy Assumptions 1, 2 and 3. Let $\calR\subset\Real^n$ be a bounded region in which we compute $V(\bfx)$. 
\begin{enumerate}
\item {\it Generating training data}. Randomly generate a set of samples $\calS_0\subset\calR$. For each $\bfx\in\calS_0$, compute $V(\bfx)$ using Algorithm 1. The outcome, $(\bfx, V(\bfx))$, is a data point. In addition, the computation results in a trajectory, $(\bfx(t), z(t))$, of the augmented system (\ref{eq_sys_aug}). One can use this trajectory to collect more data points
\EQ
\label{eq_data_midpoint}
(\bfx(t_k), V(\bfx)-z(t_k)), &1\leq k\leq K
\EE
where $\{ t_k >0\}$ is a sequence of time to be determined depending on the number and location of data points one may need along a single trajectory. The final set of data is denoted by $\calS_{training}$. 
\item {\it Generating validation data}. Repeat the computation in Step 1 to generate another set of data, $\calS_{validation}$. 
\item {\it Neural network and loss function}. Design a feedforward neural network (Figure \ref{fig_NNplot}), 
\EQ
\label{eq_NN}
V^{NN}(\bfx) = g_M\circ g_{M-1}\circ \cdots g_1(\bfx)
\EE
where $g_k(\bfz)=\sigma(\bfW_k\bfz +\bfb_k)$, $\bfz$ is a vector, $\sigma$ is a vector-valued activation function such as hyperbolic tangent, logistic or ReLU function. For the purpose of training, a widely used loss function is 
\EQ
l(\bfW,\bfb)=\Fr{1}{\abs{\calS_{training}}}\ds\sum_{\bfx\in \calS_{training}} (V^{NN}(\bfx)-V(\bfx))^2.
\EE
\item {\it Training and validation}. Find the parameters $\bfW$ and $\bfb$ for $V^{NN}(\bfx)$ that minimize $l(\bfW,\bfb)$. Using $\calS_{validation}$ to empirically evaluation the accuracy of $V^{NN}(\bfx)$. 
\end{enumerate}
There are many platforms and software packages one can use to conduct the training process for the problem of supervised learning formulated in Algorithm 2. In this paper, the training is conducted using TensorFlow \cite{abadi}.

\section{An applied example - power systems}
\subsection{Power system model}
\begin{figure}
\begin{centering}
\includegraphics[width = 2.75in]{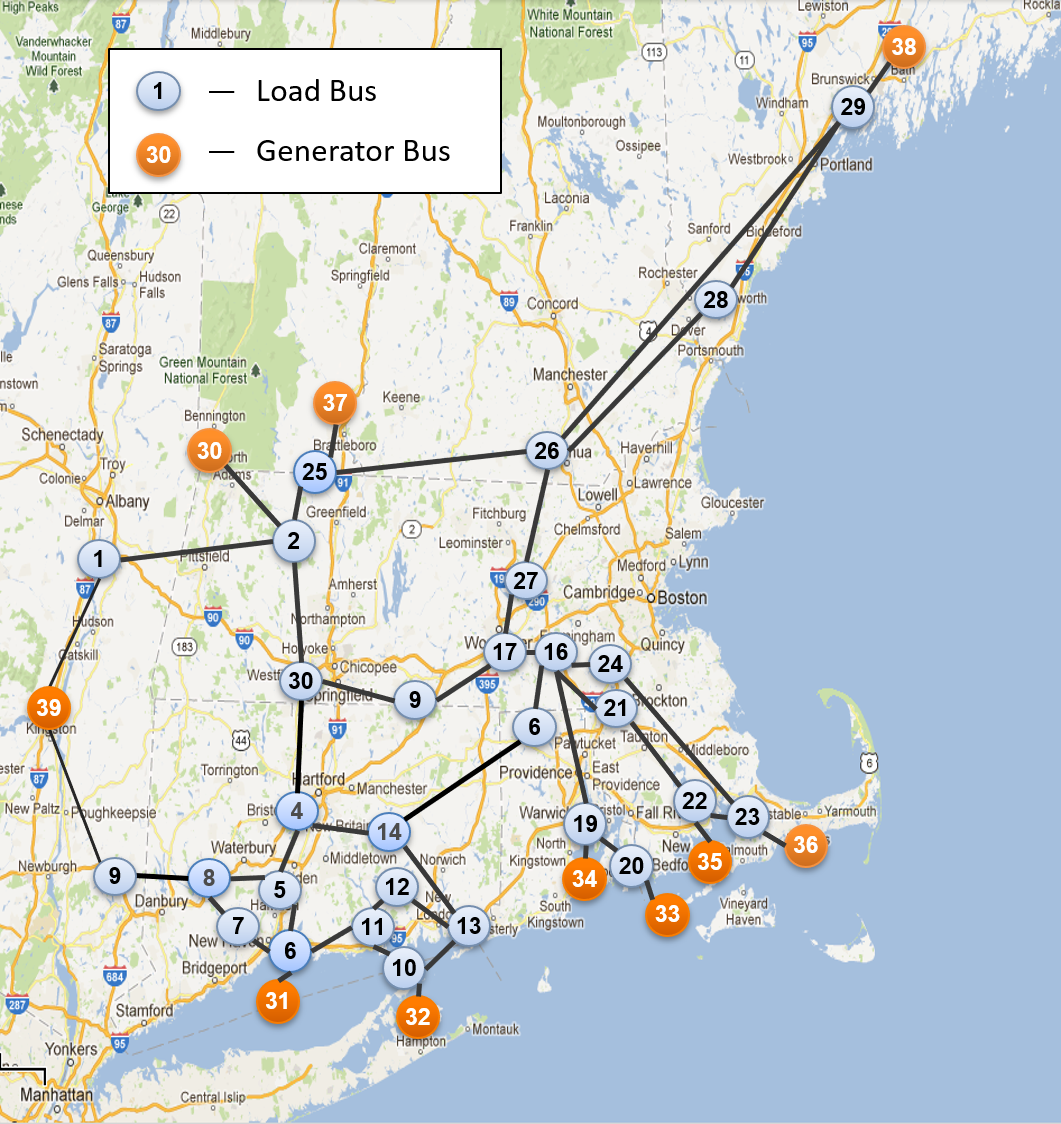}
\par\end{centering}
\caption{\label{fig:NE-39}New England 10-generator 39-bus System}
\end{figure}

In this section, we find a solution to Zubov's equation that characterizes the DOA of a 10-generator
39-bus power system.
The system's topology is shown in Figure \ref{fig:NE-39}. It is a reduction of the power grid in New England area that retains
39 main power substations and the backbone transmission
network by its 46 branches. Among the 39 substations, there are 10 generators. The first nine generators are placed at buses
No. 30 to No. 38 representing nine main power plants in the area and
the tenth generator placed at No. 39 is an equivalent of the neighboring
New York power grid. Effectively finding the stability of large scale power systems in real-time is an open problem. Illustrated in this section, the proposed Algorithm 2 provides a method of checking stability based on deep learning for which the accuracy can be empirically validated. The mathematical model of the testing system is
given by equation (\ref{eq_powersystem}), in which each of the 10 generators
is modeled by two first order ODEs, so-called swing equations. For the $i$th generator, the two state variables in  (\ref{eq_powersystem}) 
are $\delta_{i}$, the rotor angle in radian,  and $\omega_{i}$, the rotor speed
in radian per second. Other parameters include $H_{i}$ (the inertial constant of the generator),
$\omega_{0}=2\pi\times f_{0}$ (the synchronous
angular frequency in radian per second  for an ac power system with frequency $f_{0}$),
$D$ (the damping coefficient), $P_{m}$ (the mechanical power
input from the turbine), $E_{i}$ (the electromotive force or internal
voltage of the generator). In addition, $G_{ij}+jB_{ij}$,  the mutual admittance
between $E_{i}$ and $E_{j}$, is the $i^{th}$ row $j^{th}$ column
element of the admittance matrix among all electromotive forces, and
$G_{ii}$ is the conductance representing the local load seen from
$E_{i}$. Details about the model and its parameters refer to \cite{athay}.
\EQ
\label{eq_powersystem}
\Fr{d\omega_{i}}{dt}  =\Fr{\omega_{0}}{2H_{i}}\left( P_{m}-D\frac{\omega_{i}-\omega_{0}}{\omega_{0}}-E_{i}^{2}G_{ii}-\sum_{j=1,j\neq i}^{10}E_{i}E_{j}[B_{ij}\sin(\delta_{i}-\delta_{j})+G_{ij}\cos(\delta_{i}-\delta_{j})]\right)\\
\Fr{d\delta_{i}}{dt}  =\omega_{i}-\omega_{0},
\EE
where $i=1,...,10$.

\subsection{A feedforward neural network approximation of $V(\bfx)$}
\begin{figure}[!ht]
\centering
\includegraphics[width = 4.0in]{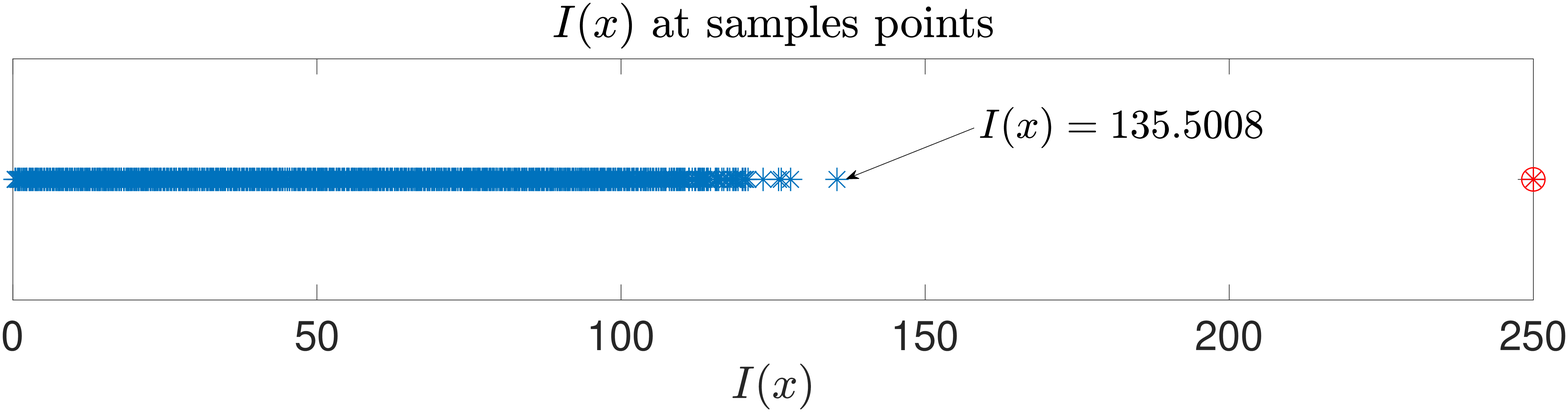}
\begin{minipage}{4.5in}
\caption{I-value plot. The cluster on the left side contains all $I(\bfx)$ satisfying (\ref{eq_dzdt}). For sample points satisfying (\ref{eq_dzdt1}), they are symbolically represented by `$ {\scriptstyle \circledast}$' at $I(\bfx)=250$.} 
\end{minipage}
\label{fig_Ix_powersystem}
\end{figure}

In the following, the state variables of (\ref{eq_powersystem}) is denoted by 
$$\bfx=\MT \omega_1 &\delta_1 &\omega_2 &\delta_2 \cdots &\omega_{10} &\delta_{10}\EM^\top\in\Real^{20}.$$ The vector field on the right-hand side in (\ref{eq_powersystem}) is denoted by $\bff(\bfx): \bfx\in \Real^{20} \rightarrow \Real^{20}$.
The computation is around a fixed equilibrium $\bfx_0$ in which
\EQ
\omega_i=120\pi\\
\delta_0 \approx  [ -0.0335, 0.0470, 0.1586, 0.1641, 0.1114, 0.1726, 0.2220, 0.1243, 0.2723,-0.1726]^\top
\EE
The sample trajectories have random initial states following the uniform distribution in $\calR$,
\EQ
\calR = \{ \bfx\in \Real^{20} ; \; -0.4\pi < \delta_i-(\delta_{0})_i<0.4\pi, 1.5<\bfw_i-120\pi<1.5 \mbox{ for } 1\leq i\leq 20\}
\EE
The system has infinitely many equilibria. In fact, the set of equilibrium, $\calE$, is a line in $\Real^{20}$. We use the following $W(\bfx)$ in the computation of $I(\bfx)$,
\EQ
W(\bfx)=\Fr{1}{1000}\bff^\top(\bfx)\bff(\bfx)
\EE
The weight in this function is for rescaling $I(\bfx)$. Otherwise, its value gets very large as the initial state approaching the boundary of $\calD$. For Algorithm 1, we choose $M=250$. In the I-value plot (Figure \ref{fig_Ix_powersystem}) based on $N=28,092$ samples,  this value of $M$ has a clear gap from the cluster of $I(\bfx)$ on the left side. Corresponding to the boundary layer $\calD_{layer}(\tanh(20))$,  the scale factor for $V(\bfx)$ should be $\alpha = 20/250=0.08$. 

\begin{figure}[!ht]
\centering
\includegraphics[width = 3.25in]{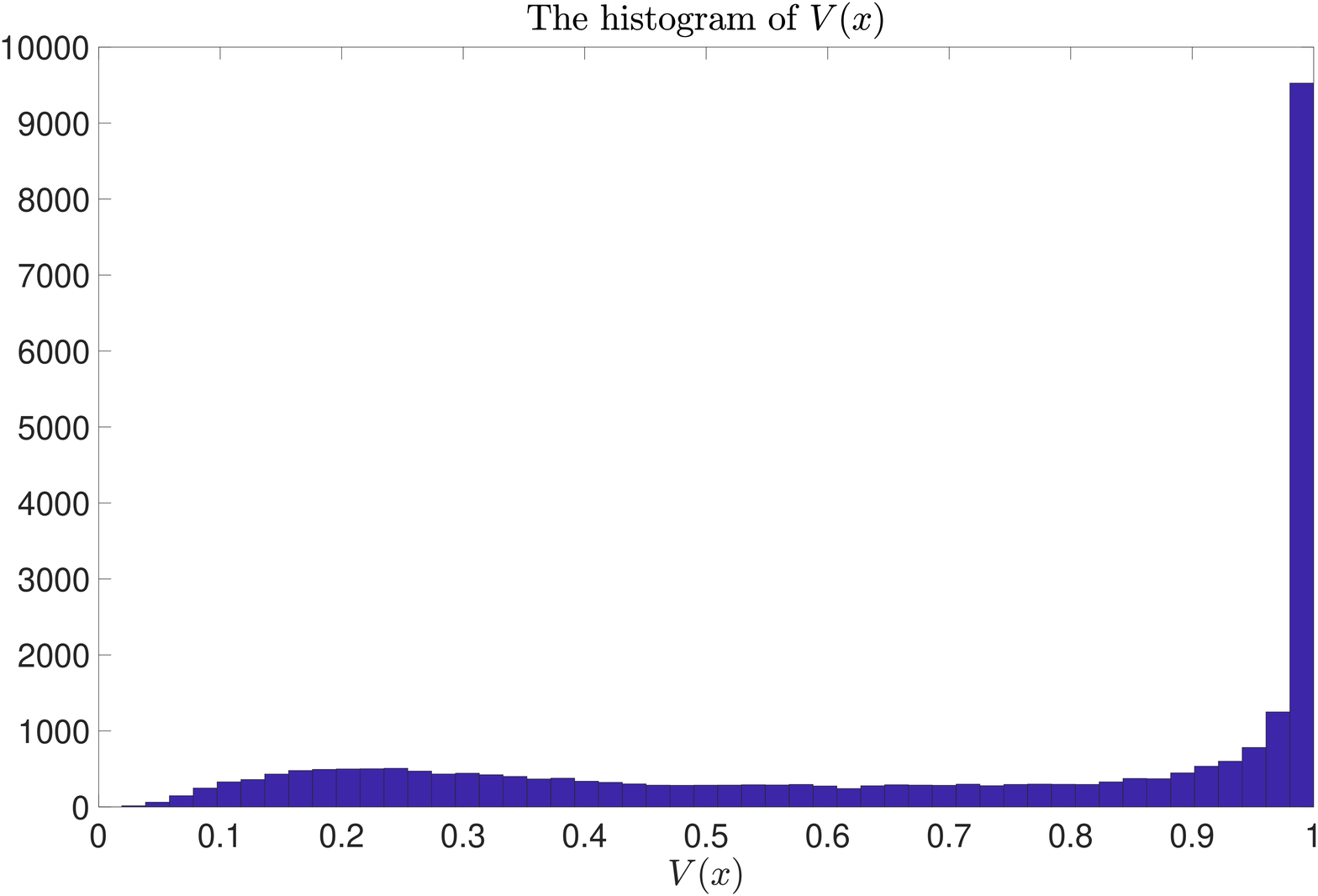}
\caption{ The histogram of $V(\bfx)$, $\bfx\in \calD_{training}$. }
\label{fig_V_hist}
\end{figure}

The data sets for training and validation are generated using Algorithm 1. For each data set, a total of $6,000$ trajectories are computed by solving the augmented system. Along each stable trajectory, four more data points are generated using (\ref{eq_data_midpoint}). The sizes of data sets are
\EQ
| \calS_{training} |= 28,092, & |\calS_{validation}|=28,100.
\EE
The level surfaces of $V(\bfx)$ form the boundaries of invariant sets in $\calD$. The data sets should contain information representing all level surfaces. The graph in Figure \ref{fig_V_hist} is the histogram of $V(\bfx)$ for $\bfx \in \calD_{training}$. It shows that the training data distributed evenly from $V(\bfx)=0.1$ to $V(\bfx)=0.9$ without missing any segment. A lot of data is concentrated near $V(\bfx)=1$. This is reasonable because $V(\bfx)=1$ represents the boundary of DOA, where we would like to provide reliable prediction. 

\begin{figure}[!ht]
\centering
\includegraphics[width = 3.25in]{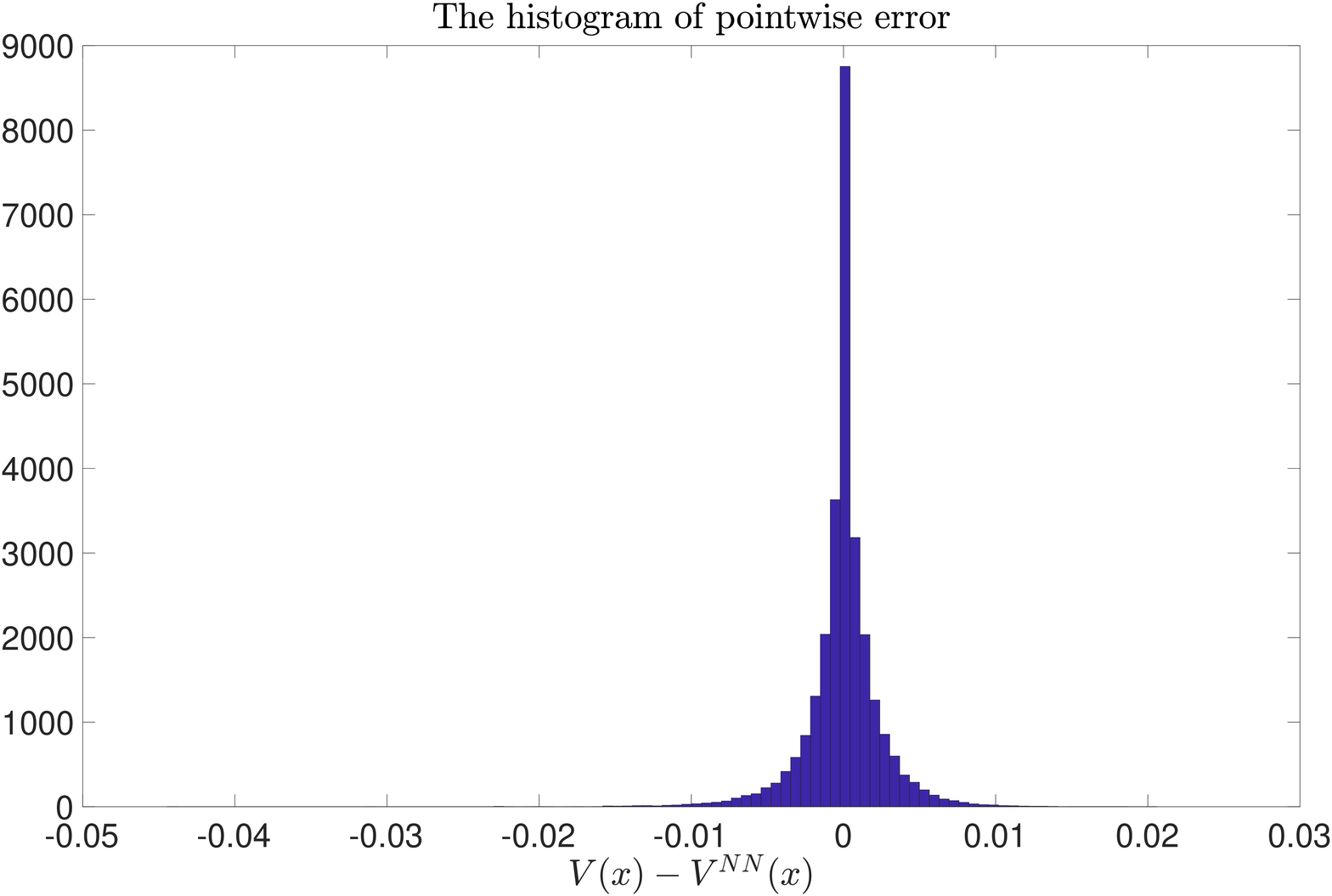} \includegraphics[width = 3.25in]{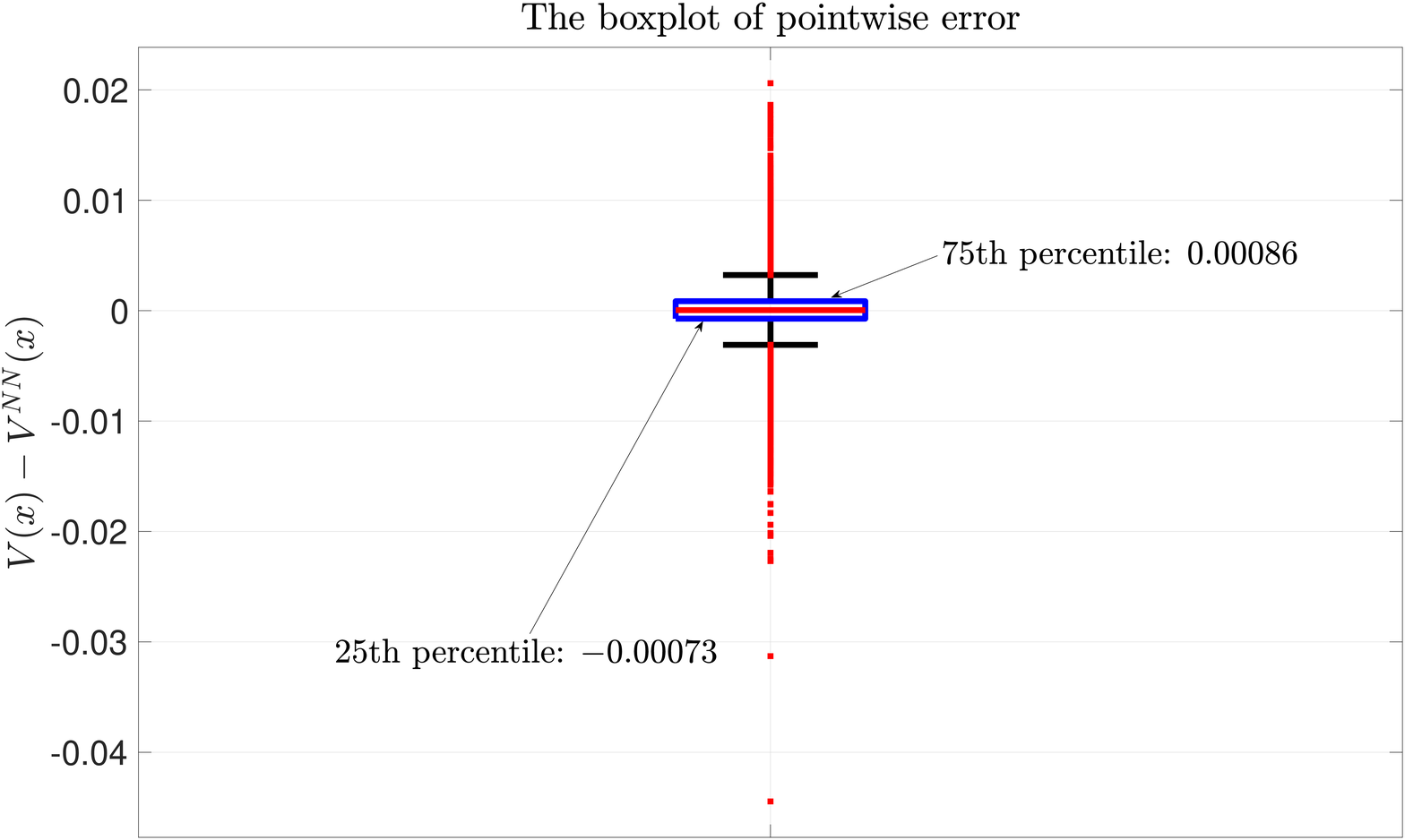}
\caption{ The histogram and boxplot  of $V(\bfx)-V^{NN}(\bfx)$, $\bfx \in \calD_{validation}$.  }
\label{fig_err_hist}
\end{figure}

Figure \ref{fig_NNplot} shows a feedforward neural network defined in (\ref{eq_NN}). In this example, the neurons in hidden layers are defined by $\sigma (x)=\tanh (x)$. The output layer is a linear function. The input of the neural network is the state of the power system, $\bfx\in \Real^{20}$. The output is $V^{NN}(\bfx)$ that approximates $V(\bfx)$. The neural network has $16$ hidden layers; each layer has $40$ neurons. The training is conducted using TensorFlow \cite{abadi}. The accuracy is measured using root-mean-square error (RMSE)
\EQ
RMSE = \sqrt{\Fr{\ds\sum_{\bfx\in \calD_{validation}}\left(V^{NN}(\bfx)-V(\bfx)\right)^2}{|\calD_{validation}|}}
\EE
The RMSE of $V^{NN}(\bfx)$ in this example is
\EQ
RMSE = 2.0\times 10^{-3}
\EE
The pointwise error over $\calD_{validation}$ has a bell shaped distribution. The histogram and boxplot of pointwise error 
$$
\{ V^{NN}(\bfx)-V(\bfx); \;\; \bfx\in \calD_{validation}\}
$$
are shown in Figure \ref{fig_err_hist}. The error variation is small. The value of the $75$th percentile is $0.00086$ and the $25$th percentile is $-0.00073$. The maximum pointwise error is $0.0445$.\\

\noindent \textbf{Remark} The augmented system (\ref{eq_sys_aug}) has multiple variables including $z$. The result of this section shows that a neural network can provide an approximation of $z$ without solving the entire system of equations (after the network is trained). This idea is promising in other applications such as data assimilation of systems that require high computational costs. In \cite{otsuka,shi}, for instance, deep learning is applied to nowcasting, i.e., predicting the future rainfall intensity in a local region over a relatively short period of time without running a full-scale data assimilation system. In general, applying deep learning in data assimilation for a local or partial approximation of large-scale systems is a new direction of research with potential applications to many areas. 

\subsection{Is it possible to break the curse of dimensionality?}
\begin{figure}[!ht]
\centering
\includegraphics[width = 3.25in]{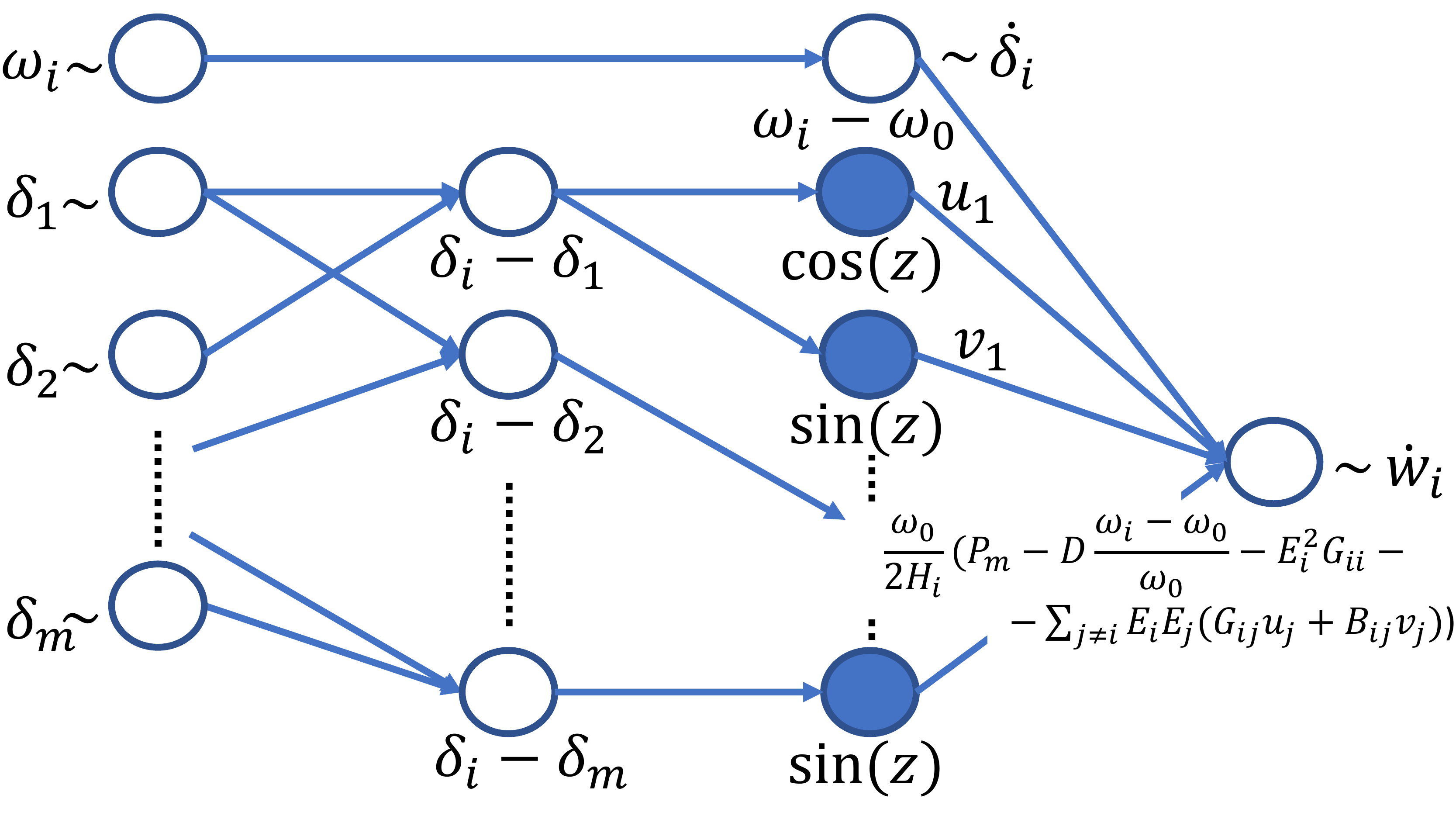} 
\caption{ The layered DAG of the function in (\ref{eq_powersystem}) as a compositional function}
\label{fig_Pei}
\end{figure}

\begin{figure}[!ht]
\centering
\includegraphics[width = 3.25in]{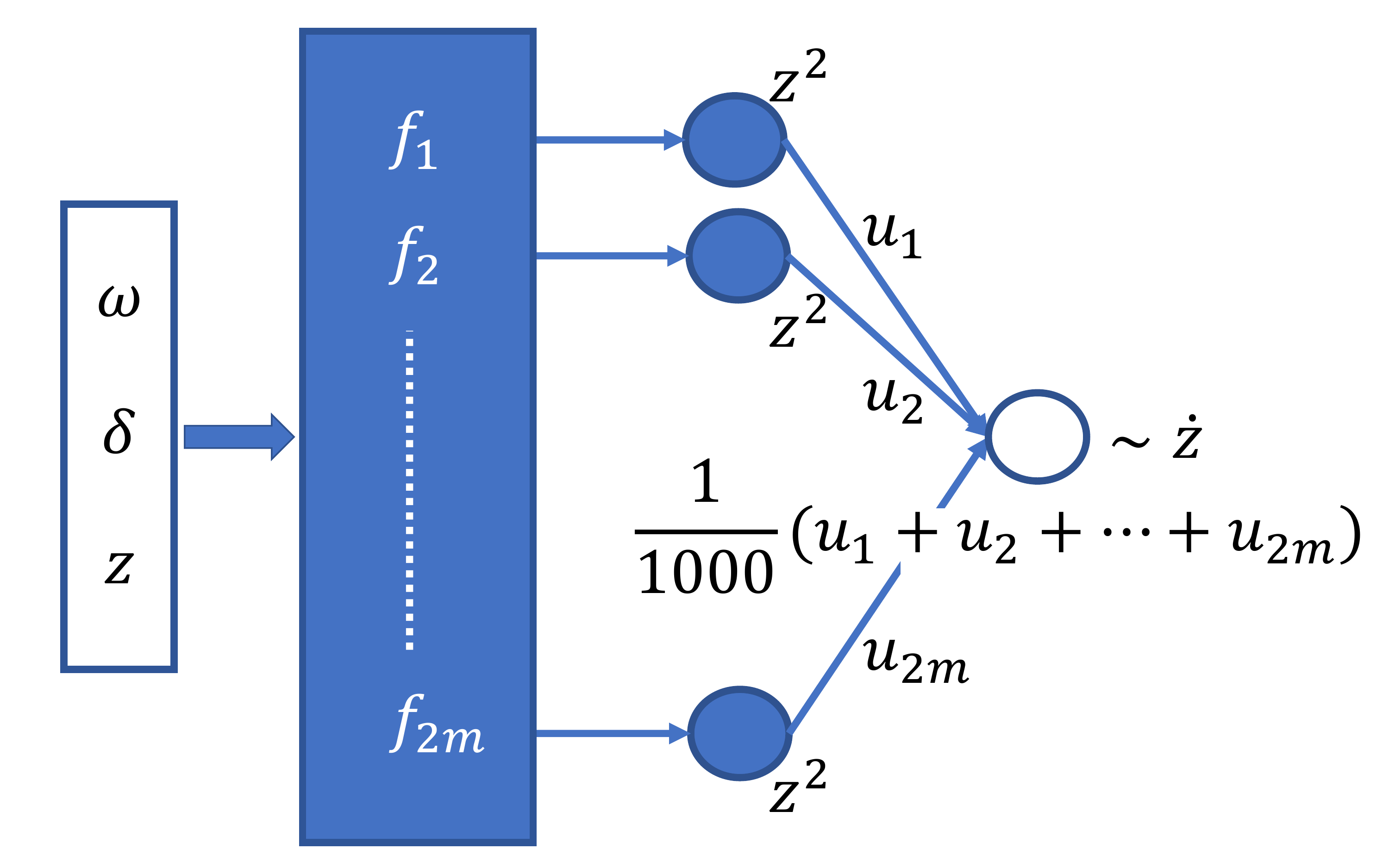} 
\caption{ The layered DAG of the function in (\ref{eq_augmented}) as a compositional function}
\label{fig_z}
\end{figure}

As a solution of a partial differential equation, $V(\bfx)$ is very difficulty to compute, if not impossible, when the state space has a high dimension. For example, a PDE discretization based on $N$ points in each dimension requires a grid that has $N^n$ points in a $n$-dimensional state space. For the power system (\ref{eq_powersystem}) in which $n=20$, the grid size is $10^{40}$ if $N=100$. Any numerical method based on a grid of such size is practically intractable. Other numerical methods of solving PDEs, such as the finite element method, have the same problem if the state space has a high dimension. The fundamental challenge lies in the fact that the complexity of numerical approximation increases with the dimension of state space exponentially. This phenomenon is called the {\it curse of dimensionality}.  In \cite{kanggong}, it is proved that the approximation of compositional functions using neural networks may not suffer the curse of dimensionality. Rather than an exponential function, the approximation error depends on the dimension as a polynomial. The model in (\ref{eq_powersystem}) is defined based on compositional functions, in the same form as those addressed in \cite{kanggong}. Its compositional structure can be represented using a layered directed acyclic graph (layered DAG), which is shown in Figure \ref{fig_Pei}. The augmented system has an additional equation
\EQ
\label{eq_augmented}
\dot z= W(\bfx)=\Fr{1}{1000}\bff^\top(\bfx) \bff(\bfx), & z(0)=0,
\EE
where $\bff$ is the vector field representing the right-hand side in (\ref{eq_powersystem}). The function in (\ref{eq_augmented}) is also a compositional function, whose layered DAG is shown in Figure \ref{fig_z}.

Applying the theory in \cite{kanggong}, the colored nodes in Figures \ref{fig_Pei} and \ref{fig_z} represent nonlinear functions. The white nodes represent linear functions. The compositional features, defined for nonlinear nodes only,  are briefly introduced as follows. 
\begin{itemize}
\item  $| \calV |$ (complexity feature): The total number of nonlinear nodes in the layered DAG, where $\calV$ is the set of nonlinear nodes of the compositional function.
\item $ r_{max}$ (dimension feature): In this paper we assume that all nodes in $\calV$ are $C^1$ functions. The dimension feature is defined to be the largest input dimension of nodes in $\calV$.  For the compositional function in Figure \ref{fig_Pei} and \ref{fig_z}, $r_{max}=1$. 
\item  $\Lambda$ (volume feature): Each nonlinear node, as a function, has a domain. Assume that the domain is a square of edge length $R$. The volume feature is defined to be the largest value in 
\EQ
\{ \max\{R,1\} || f ||; \;f\in\calV\}, 
\EE
where $|| \cdot ||$ is the Sobolev norm 
\EQ
\label{eq_sobolevnorm}
|| f || = ||f||_{L^\infty}+\ds\sum \norm{\Fr{\partial f}{\partial x_i} }_{L^\infty}.
\EE 
\item $L_{\max}$ (Lipschitz constant feature): It is the largest Lipschitz constants associated with nonlinear nodes. Note that this Lipschitz constant is defined based on the layered structure of the DAG. Discussions in the rest of this session have some further explanations. For more details, the readers are referred to \cite{kanggong}. 
\end{itemize}
For the augmented system (\ref{eq_powersystem})-(\ref{eq_augmented}), let us denote the augmented vector field by $\bar \bff = [\bff(\bfx)^\top \; W(\bfx)]^\top$. Let $\bar \bfphi(t,\bfx)$ represents the solution of (\ref{eq_powersystem})-(\ref{eq_augmented}) in which $\bfx$ is the initiate state.  Suppose that $\bar \bfphi(t,\bfx)$, for $t\in [0, T]$, is contained in a bounded closed set in the state space. Given a time interval $[0, T]$, it is proved in \cite{kanggong} that there always exists a deep feedforward neural network, denoted by $\bar \bfphi^{NN}(\bfx)$, in which  activation functions are $C^\infty$. Furthermore, the network satisfies
\EQ
\label{eq_NN_error}
\norm{\bar \bfphi^{NN}( \bfx)-\bar \bfphi(T, \bfx)}_2 < (C_1 L_{max}\Lambda \abs{\calV}+C_2)\tilde n^{-1/r_{max}}
\EE
where $\tilde n$ is an integer that determines the total number of neurons in $\bar\bfphi^{NN}(\bfx)$, i.e., the complexity of the neural network,
\EQ
\label{eq_NN_complexity}
\mbox{ The complexity  of } \bar\bfphi^{NN} \leq \left(\tilde n^{1/r_{max}}+1\right)\tilde n\abs{\calV}
\EE
The constants, $C_1$ and $C_2$, in (\ref{eq_NN_error}) are determined by $\norm{\bar\bff}_2$, $\norm{\Fr{\partial \bar\bff}{\partial \bar \bfx}}_2$, $T$ and the input dimensions of the nodes in $\calV$. 

From Algorithm 1, an approximate solution to Zubov's equation can be represented in the following form
\EQ
\label{eq_Vx_T}
V(\bfx)=\tanh(\alpha z(T,\bfx)),
\EE
where $z(t,\bfx)$ is the solution of the augmented system (\ref{eq_sys_aug}) and $T>0$ is a constant such that $z(t,\bfx)$ satisfies the criteria (\ref{eq_dzdt1}) or (\ref{eq_dzdt}) for $\bfx$ in a bounded set. 

\begin{theorem}
\label{thm_compositional}
Consider a power system  (\ref{eq_powersystem}) that has $m$ generators and its Lyapunov function $V(\bfx)$ defined in (\ref{eq_Vx_T}). Let $\calR\subset \Real^{2m}$ be a bounded set. Let $V^{NN}( \bfx)$ represent a feedforward neural network that has $n^{NN}$ neurons, which are hyperbolic tangent functions. Then, there exists $V^{NN}( \bfx)$ that satisfies
\EQ
\label{eq_NN_err_power}
\abs{V^{NN}( \bfx)-V(\bfx)} < (C_1 m^2+C_2)\Fr{m}{\sqrt{n^{NN}}}
\EE
for $\bfx\in\calR$, where $C_1$ and $C_2$ are constants independent of $m$. They depend on $\norm{\bar\bff}_2$, $\norm{\Fr{\partial \bar\bff}{\partial \bar \bfx}}_2$, $T$, $\alpha$ and the parameters in (\ref{eq_powersystem}). 
\end{theorem}

\noindent{\it Proof}. The proof is based on (\ref{eq_NN_error}) - (\ref{eq_NN_complexity}). We first prove some properties of compositional features. Specifically, $L_{max}$ and $\Lambda$ are independent of $m$; and $\abs{\calV}$ is a polynomial function of $m$. Following the definition in \cite{kanggong}, the Lipschitz constant associated with a node is the Lipschitz constant of the function with respect to the node when the node is treated as a variable. For example, consider the function in Figure \ref{fig_Pei}, in which $\cos(z)$ and $\sin(z)$ are the nonlinear nodes. If one of them, for instance the first $\cos (z)$ connecting to $\delta_i-\delta_1$,  is treated as a variable, the Lipschitz constant of the function associated with this node equals 
\EQ
\Fr{\omega_0}{2H_i}E_iE_1G_{i1}.
\EE 
The Lipschitz constant of nonlinear nodes in Figure \ref{fig_z} equals $1/1000$. Therefore, $L_{max}$ is independent of $m$. The value of $\Lambda$ depends on the radius of the domain of nonlinear nodes (the sine and cosine functions in (\ref{eq_powersystem})) and their Sobolev norm (\ref{eq_sobolevnorm}). They are independent of $m$. The set of nonlinear nodes, $\calV$, consists of $\cos(z)$, $\sin(z)$, and $z^2$ for the functions in (\ref{eq_NN_error}) - (\ref{eq_NN_complexity}). For each $1\leq i\leq m$, there are $2(m-1)$ nonlinear nodes for the function in Figure \ref{fig_Pei}. There are $2m$ nonlinear nodes for the function in Figure \ref{fig_z}.  Therefore, the total number of nonlinear nodes in the function of the augmented system (\ref{eq_powersystem})-(\ref{eq_augmented}) is $\abs{\calV}=2m^2$. The dimension feature is $r_{max}=1$ because all nonlinear nodes have a single input. Let $\bar \bfphi (t, \bfx)$ represent the solution of the augmented system (\ref{eq_powersystem})-(\ref{eq_augmented}). The inequalities (\ref{eq_NN_error}) and (\ref{eq_NN_complexity}) for the power system are equivalent to
\EQ
\label{eq_err_tmp}
\norm{\bar \bfphi^{NN}( \bfx)-\bar \bfphi(T, \bfx)}_2 < (\bar C_1m^2+\bar C_1) \Fr{1}{\tilde n}\\
n^{NN} < 2\tilde n^2(2m^2)
\EE
for some feedforward neural network, $\bar \bfphi^{NN}(\bfx)$, whose complexity is $n^{NN}$. Solving for $\tilde n$ from the second equation in (\ref{eq_err_tmp}) and then applying it to the first equation, we have
\EQ
\label{eq_err_tmp1}
\norm{\bar \bfphi^{NN}( \bfx)-\bar \bfphi(T, \bfx)}_2 < (\bar C_1m^2+\bar C_1) \Fr{2m}{\sqrt{n^{NN}}}\\
\EE
Define 
$$V^{NN}(\bfx)=\tanh (\alpha \bar\bfphi^{NN}_{2m+1}(\bfx)).$$
where $\bar\bfphi^{NN}_{2m+1}(\bfx)$ is the last component of $\bar\bfphi^{NN}(\bfx)$. Note that $z(t,\bfx)$ is the last component of $\bar\phi(\bfx)$. Because the composition of neural networks results in a neural network and because $\tanh(\cdot)$ is Lipschitz, it is straightforward to justify that (\ref{eq_err_tmp1}) implies (\ref{eq_NN_err_power}). 
$\blacksquare$\\

Theorem \ref{thm_compositional} implies that the neural network approximation of $V(\bfx)$ has an error that is a cubic polynomial of  $n=2m$, the state space dimension. The error does not grow exponentially. We can conclude that there exists a neural network approximate solution to Zubov's equation that breaks the curse of dimensionality. 

We would like to emphasize that finding the desired neural network in Theorem \ref{thm_compositional} is an open problem. In neural network training and validation, no matter how large is the data set, there is no general methodologies that can guarantee the error upper bound (\ref{eq_NN_err_power}) for all $\bfx\in \calR$. In addition, the training process is a nonconvex optimization. It may converge to a local minimum; and the resulting neural network may not provide an accurate approximation of $V(\bfx)$. These drawbacks, however, are challenges facing general practice of deep learning, not limited only to the problems studied in this paper. 

\section{Conclusions}
Characterizing the DOA for general nonlinear systems of ODEs has been an open problem for decades. In this paper, we develop practical methods of describing the DOA by finding a Lyapunov function that solves Zubov's equation. The theoretical foundation is a theorem in which we construct and prove a solution to Zubov's equation. This solution can be represented in an integral form that is easy to compute. All that is needed is an ODE solver. The parameters in the algorithm can be determined based on simulation data. For problems that ODE solvers do not meet the requirement of real-time computation, such as large scale power grids, this algorithm can still be used to generate data so that a feedforward neural network is trained to approximate the Lyapunov function. In fact, we mathematically prove that the swing equation of power systems admits a neural network approximate solution to Zubov's equation. The approximation error and complexity of the neural network are polynomial functions of the number of generators. It implies that a solution exists that breaks the curse of dimensionality. The deep learning approach is applied to the New England 10-generator 39-bus power system, an interesting research problem in its own right.

\end{document}